\documentclass{amsart}
\title{Almost all strongly quasipositive braid closures are fibered}
\author{Ian Banfield}
\address{Department of Mathematics, Boston College, Chestnut Hill, MA-02467}
\email{ian.banfield@bc.edu}
\usepackage{graphicx}
\usepackage{amsfonts}
\usepackage{amsmath}
\usepackage{pstricks}
\usepackage{todonotes}
\usepackage{subcaption}
\usepackage{amssymb}
\usepackage{amsthm}
\usepackage{enumitem}
\usepackage{ tipa }

\newcommand{\ul}[1]{\underline{#1}}

 \newcommand{\executeiffilenewer}[3]{%
 	\ifnum\pdfstrcmp{\pdffilemoddate{#1}}%
 	{\pdffilemoddate{#2}}>0%
 	{\immediate\write18{#3}}\fi%
 }
 \IfFileExists{/dev/null}{%
 	\newcommand{\Inkscape}{inkscape }%
 }{
 \newcommand{\Inkscape}{"C:/Program Files (x86)/Inkscape/inkscape.exe" }%
}

\theoremstyle{plain} \numberwithin{equation}{section}
\newcounter{dummy} 
\numberwithin{dummy}{section}
\newtheorem{theorem}[dummy]{Theorem}
\newtheorem{lemma}[dummy]{Lemma}
\newtheorem{cor}[dummy]{Corollary}

\newtheorem{fact}[dummy]{Fact}
\newtheorem{definition}[dummy]{Definition}

\makeatletter
\newtheorem*{rep@theorem}{\rep@title}
\newcommand{\newreptheorem}[2]{%
	\newenvironment{rep#1}[1]{%
		\def\rep@title{#2 \ref{##1}}%
		\begin{rep@theorem}}%
		{\end{rep@theorem}}}
\makeatother
\newreptheorem{theorem}{Theorem}
\newreptheorem{lemma}{Lemma}

\begin{document}
\newpage
\begin{abstract}
	We use the Birman-Ko-Lee presentation of the braid group to show that all closures of strongly quasipositive braids whose normal form contains a positive power of the dual Garside element $\delta$ are fibered. We classify links which admit such a braid representative in geometric terms as boundaries of plumbings of positive Hopf bands to a disk. Rudolph constructed fibered strongly quasipositive links as closures of positive words on certain generating sets of $B_n$ and we prove that Rudolph's condition is equivalent to ours. Finally, we show that the braid index is a strict upper bound for the number of crossing changes required to fiber a strongly quasipositive braid.
\end{abstract}

\maketitle



\section{Introduction}\label{introduction}

\begin{figure}
			\begin{subfigure}[t]{.6\textwidth}
				\centering
	\executeiffilenewer{SQPbraidedWithSeifertSurface.svg}{SQPbraidedWithSeifertSurface.pdf}{%
		\Inkscape -z -D --file="SQPbraidedWithSeifertSurface.svg" --export-pdf="SQPbraidedWithSeifertSurface.pdf" --export-latex}%
	\scalebox{.4}{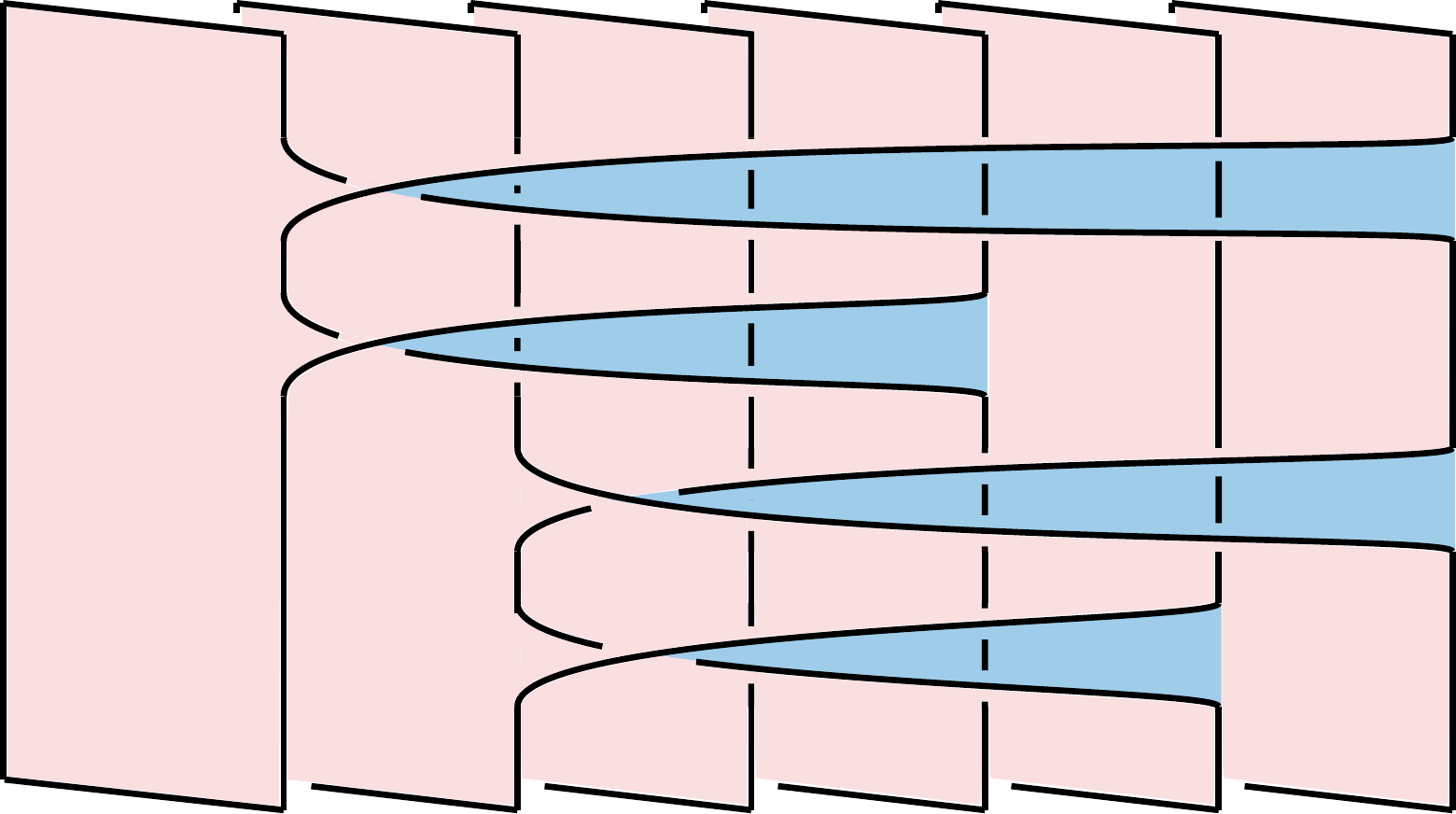}%

				\subcaption{Bennequin surface $\Sigma_w$}
				\label{fig:sqpbraid}
			\end{subfigure}
			\begin{subfigure}[t]{.3\textwidth}
				\centering
	\executeiffilenewer{ChargedFenceDiagram.svg}{ChargedFenceDiagram.pdf}{%
		\Inkscape -z -D --file="ChargedFenceDiagram.svg" --export-pdf="ChargedFenceDiagram.pdf" --export-latex}%
	\scalebox{.4}{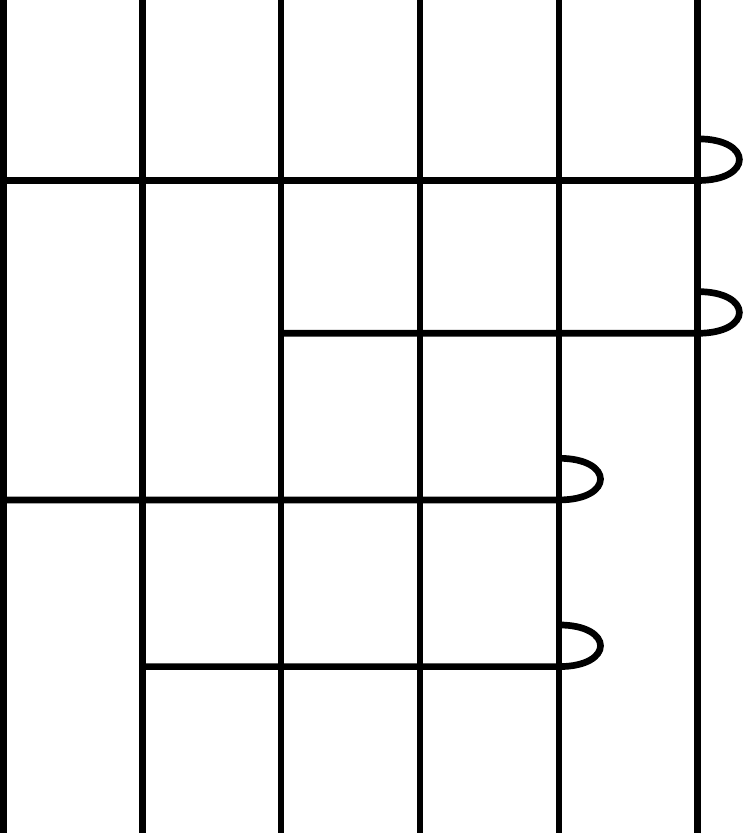}%

				\subcaption{Charged fence diagram}
				\label{fig:chargedfencediagram}
			\end{subfigure}
			\caption{The closure of the strongly quasipositive braid $\beta = a_{6,1} a_{6,3} a_{5,1} a_{5,2}$ as the boundary of a quasipositive surface.}
		\end{figure}

Strongly quasipositive links, i.e. links that are representable as closures of strongly quasipositive braids, are an interesting class of links first studied by Rudolph, cf.  \cite{MR1452826}. Geometrically, they can be described as boundaries of quasipositive surfaces in $S^3$, see Figure \ref{fig:sqpbraid}.

It is known that the slice Bennequin inequality \cite{MR1193540} is sharp for quasipositive links in general; indeed for strongly quasipositive links the Bennequin surface is actually embedded in $S^3$, which implies that the slice genus equals the genus of the link \cite{MR1241873,MR1193540}. This can be used to combinatorially compute the genus of a strongly quasipositive link from any strongly-quasipositive braid representation. It also implies that the only slice strongly quasipositive knot is the unknot. There are connections to invariants coming from knot homology theories; e.g. the slice genus bounds from Rasmussen's $s$-invariant in Khovanov homology and the $\tau$-invariant in knot Floer homology are sharp \cite{MR2384833}.

Particularly interesting are fibered strongly quasipositive links: By work of Hedden and Rudolph, it is known that these are exactly the fibered links that induce the unique tight contact structure on $S^3$ \cite{MR2646650}. Knot Floer homology detects if a link is fibered \cite{MR2357503}, and further if a fibered link can be represented as a strongly quasipositive braid closure \cite{MR2646650}.

We will show that ``almost all" strongly quasipositive braid closures are fibered.  Our proof uses a presentation of the braid group due to Birman-Ko-Lee, associated to the so-called dual Garside structure on the braid group. 

\begin{replemma}{thm:mainlemma}
	Let $\beta = \delta P$, where $\delta = \sigma_{n-1} \sigma_{n-2} \dots \sigma_{1}$ is the Dual Garside element and $P$ is a BKL-positive word. Then the braid closure $\hat{\beta}$ is fibered.
\end{replemma}

See Section \ref{braids} for the definition of a BKL-positive word and the Dual Garside normal form. We prove the following theorem.

\begin{reptheorem}{thm:maintheorem}
	A link $L$ which can be represented as the closure of a braid whose normal form contains a positive power of the Dual Garside element is fibered.
\end{reptheorem}

In Section \ref{section:hopfplumbedbaskets} we classify which fibered strongly quasipositive braids arise from our construction and relate the condition that a link $L$ is represented as a braid closure whose normal form contains the Dual Garside element to a construction of fibered links due to Rudolph \cite{MR1857666}.

\begin{reptheorem}{thm:hopfplumbedbaskets}
	A link $L$ is the boundary of a plumbing of positive Hopf bands to a disk $D$ along arcs $\alpha_i \subset D$ if and only if $L$ admits a strongly quasipositive representative $\beta \in B_n$ which contains the Dual Garside element $\delta$.
\end{reptheorem}

In particular, we show that positive braid closures are of this form.

\begin{reptheorem}{thm:corpositivebraids}
	A non-split positive braid link $L$ is the closure of a strongly quasipositive braid $\beta \in B_n$ whose normal form contains a positive power of the Dual Garside element.
\end{reptheorem}

\subsection*{Acknowledgements}
The author is deeply grateful to his advisor Eli Grigsby for her inspiration and support, and her thoughtful suggestions and careful reading of the many drafts of this paper.
The author thanks Sebastian Baader, Peter Feller, Tao Li and Filep Misev for helpful conversations and their encouragement. Section \ref{section:hopfplumbedbaskets} was inspired by a conversation and subsequent email exchange with Lee Rudolph after meeting him at the Boston Topology Graduate Student conference at MIT in October 2016. Lastly, the author would like to thank Peter Feller for pointing out a mistake in the original proof of Corollary \ref{lemma:makingsqpfibredbyaddingcrossings}.

\section{The Dual Garside structure of the braid group}
\label{braids}

\begin{figure}
		 \begin{subfigure}[t]{.33\textwidth}
			\centering
	\executeiffilenewer{BKL-generator.svg}{BKL-generator.pdf}{%
		\Inkscape -z -D --file="BKL-generator.svg" --export-pdf="BKL-generator.pdf" --export-latex}%
	\scalebox{1}{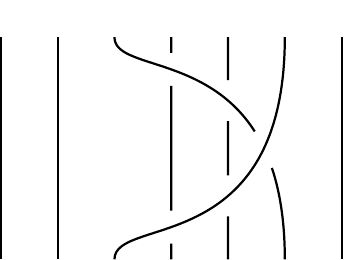}%

			\subcaption{Braid presentation}
			\label{fig:bklgenerator_braid}
		\end{subfigure}
		\begin{subfigure}[t]{.33\textwidth}
			\centering
	\executeiffilenewer{BKL-circle-presentation.svg}{BKL-circle-presentation.pdf}{%
		\Inkscape -z -D --file="BKL-circle-presentation.svg" --export-pdf="BKL-circle-presentation.pdf" --export-latex}%
	\scalebox{1}{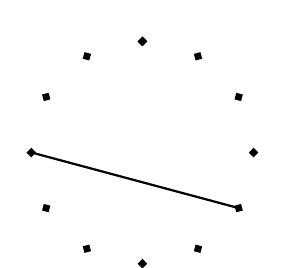}%

			\subcaption{Circle presentation}
			\label{fig:bklgenerator_circle}
		\end{subfigure}
		\begin{subfigure}[t]{.3\textwidth}
			\centering
	\executeiffilenewer{BKL-circle-presentation-example-of-braid.svg}{BKL-circle-presentation-example-of-braid.pdf}{%
		\Inkscape -z -D --file="BKL-circle-presentation-example-of-braid.svg" --export-pdf="BKL-circle-presentation-example-of-braid.pdf" --export-latex}%
	\scalebox{1}{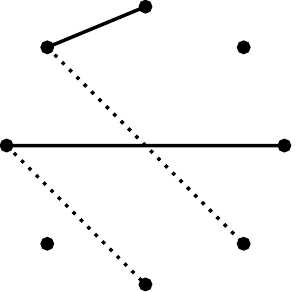}%

			\subcaption{$\beta = a_{8,1} a_{8,4}^{-1} a_{7,3} a^{-1
				}_{7,5}$}
			\label{fig:bklexample_circle}
		\end{subfigure}
		\caption{The generators $a_{i,j}$. For convenience, define $a_{i,j} = a_{j,i}$. The underlined number next to a line indicates the position in the braid word.}
		\label{fig:bklgenerator}
	\end{figure}

In \cite{MR1654165}, Birman, Ko and Lee gave a solution to the word-problem in the braid group using a new presentation of $B_n$. The generators $a_{i,j}$ in this presentation correspond to pairs of strands (see Figure \ref{fig:bklgenerator}), called band generators or Birman-Ko-Lee generators, and the relations correspond to Reidemeister moves of type II and III. 

There is a nice pictorial way to describe the generators and relators, using dots labeled $1, 2, \dots, n$ arranged in a circle \cite{MR3352238}. The generator $a_{i,j}$ is represented by a straight line between the dots labeled $i,j$, as pictured in Figure \ref{fig:bklgenerator_circle}, and its inverse by a dashed line. An element of $B_n$ is represented by a sequence of such lines, see Figure \ref{fig:bklexample_circle}. The Birman-Ko-Lee relations can be stated as follows.

\begin{figure}
		\centering
	\executeiffilenewer{BKL-circle-presentation-relations.svg}{BKL-circle-presentation-relations.pdf}{%
		\Inkscape -z -D --file="BKL-circle-presentation-relations.svg" --export-pdf="BKL-circle-presentation-relations.pdf" --export-latex}%
	\scalebox{1}{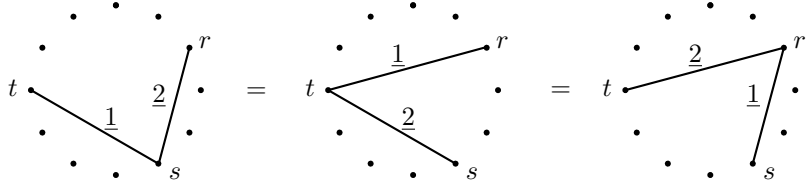}%

		\caption{Going counter-clockwise, any two consecutive boundaries of the triangle with vertices at $r,s,t$ are equivalent in $B_n$, i.e. $a_{t,s} a_{s,r} = a_{t,r} a_{t,s} = a_{s,r} a_{t,r}$.}
		\label{fig:sqprelations}
	\end{figure}

\begin{itemize}
	\item Consecutive lines in the sequence representing the braid word which do not intersect can be re-ordered arbitrarily.
	\item (``Cup Product") Consider a triangle with vertices at the dots labeled $r$, $s$, $t$. Moving counter-clockwise, any two consecutive boundaries of the triangle are equivalent. See Figure \ref{fig:sqprelations}.
\end{itemize}

\begin{figure}
		\centering
	\executeiffilenewer{Canonical-factor-from-polygon.svg}{Canonical-factor-from-polygon.pdf}{%
		\Inkscape -z -D --file="Canonical-factor-from-polygon.svg" --export-pdf="Canonical-factor-from-polygon.pdf" --export-latex}%
	\scalebox{1}{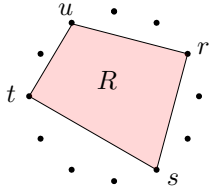}%

		\caption{The representations of the braid $\beta_R$ associated to the $4$-gon $R$ arising from cyclic counterclockwise orders of its vertices. \\ $\beta_R = a_{u,r} a_{t,u} a_{t,s} = a_{u,t} a_{t,s} a_{s,r} = a_{t,s} a_{s,r} a_{u,r} = a_{s,r} a_{u,r} a_{t,u}$.}
		\label{fig:braidfrompolygon}
	\end{figure}

The second relation implies that one can associate to a polygon $R$ a well-defined element $\beta_R \in B_n$ as follows: Let the vertices of $R$ be at $(q_1, \dots, q_r)$ listed in any cyclic counterclockwise order and let $\beta_R = a_{q_r, q_{r-1}} \dots a_{q_2,q_1}$. See Figure \ref{fig:braidfrompolygon} for an example when $R$ is a $4$-gon.

A braid $\beta$ is \textbf{BKL-positive} if it can be written as a word in only positive powers of the generators $a_{i,j}$. BKL-positive braids form a monoid, denoted by $B_n^+$, and Birman-Ko-Lee proved that two BKL-positive elements are equivalent in $B_n$ if and only if they are equivalent in $B_n^+$.

\begin{definition}
	The Dual Garside element is $\delta = a_{n,n-1} a_{n-1, n-2} \dots a_{2,1}$.
\end{definition}
As the name implies, $\delta$ is in fact a Garside element for a Garside structure on $B_n$. In the circle presentation, $\delta$ corresponds to the $n$-gon spanned by all vertices.

\begin{definition}
	The starting set $S(\beta)$ of a BKL-positive braid word $\beta$ is the set of generators $a_{i,j}$ that can appear at the start of a BKL-positive word representing $\beta$. Similarly, the finishing set $F(\beta)$ is the set of generators $a_{i,j}$ that can appear at the end of $\beta$.
\end{definition}

We will need the following facts from \cite{MR1654165}.
\begin{fact}
	The finishing set of the Dual Garside Element $\delta$ is the set of all generators, i.e. $F(\delta) = \{a_{j,i} ~|~ i, j \in \{1, \dots n\}, i \neq j\}$.
	\label{fact:finishingsetofdualgarside}
\end{fact}

\begin{fact}[Dual Garside Normal Form, \cite{MR1654165}, Theorem 3.10]
	\label{fact:bklnormalform}
	Every element $\beta \in B_n$ can be uniquely written as $\beta = \delta^k A_1 \dots A_m$, where $\delta$ is the Dual Garside element, all $A_i$ are BKL-positive and $A_i s$ is not a canonical factor for any $s \in S(A_{i+1})$.
\end{fact}

The $A_i$ in the previous theorem are called canonical factors and correspond to disjoint union of polygons in the circle presentation. The condition that $A_i s$ is not a canonical factor for any $s \in S(A_{i+1})$ ensures uniqueness of the normal form and is denoted by $A_i \lceil  A_{i+1}$.

A braid is strongly quasipositive if and only if the power of the Dual Garside element $\delta$ in the normal form is non-negative. A link is strongly quasipositive if it admits a braid representative with a non-negative power of $\delta$.

\section{Detecting fibered braids}
\label{fibrations}

A link $L$ is fibered with fiber $F$ if $\partial F = L$ and its exterior $S^3 \setminus \nu(L)$ fibers over the circle such that $F$ is a fiber. This is a strong condition; it implies that all fiber surfaces are isotopic and minimal genus. In fact, for Seifert surfaces $F$ of a fibered link the following are equivalent: $F$ is a fiber surface, $F$ is genus-minimizing and $F$ is incompressible \cite{kawauchi1996survey}.

In \cite{MR870705}, Gabai established an algorithm to detect if a link $L \subset S^3$ is fibered. The idea is as follows: Let $F$ be a minimal genus Seifert surface for $L$. The link $L$ is fibered if and only if the complementary sutured manifold $(Y, \gamma) = (S^3 \setminus (F \times I), \partial F)$ associated to $F$ is a product sutured manifold, that is, $(Y, \gamma) \simeq (\Sigma \times I, \partial \Sigma)$ for a surface $\Sigma$. Gabai proves the following theorem \cite{MR870705}.

\begin{fact}[\cite{MR870705}, Theorem 1.9]
	\label{fact:gabaifibredlinks}
	The link $L \subset S^3$ is fibered with fiber surface $F$ if and only if there exists a sequence of product disk decompositions of $(Y, \gamma) = (S^3 \setminus (F \times I), \partial F)$ that terminates in the trivial sutured manifold $(B^3, \alpha)$, where $\alpha$ is a single curve on $\partial B^3$.
\end{fact}

We will apply this theorem to a braid closure $\hat{\beta}$ of a braid word expressed in the Birman-Ko-Lee generators. Recall from Section \ref{braids} that we can express a braid $\beta = w$ as a word $w$ in the Birman-Ko-Lee generators. For such a word $w$, we can construct a canonical Seifert surface, $\Sigma_w$, called the Bennequin surface, as follows.

\begin{itemize}
	\item Draw $n$ parallel disks, one for every strand of the braid $\beta$.
	\item For each generator $a_{j,i}$ (inverse of $a_{j,i}$, resp.) attach a positively (negatively, resp.) twisted band between the disks corresponding to strands $i$ and $j$, going over all other disks.
\end{itemize}

See Figure \ref{fig:sqpbraid} for an example of the Bennequin surface for $\beta = a_{6,1} a_{6,3} a_{5,1} a_{5,2}$. For strongly quasipositive braid closures, the Bennequin surface is minimal genus among all smoothly embedded surfaces in $B^4$ bounded by the oriented link; this follows from Rudolph's proof of the slice Bennequin inequality \cite{MR1193540}, which in turn relies on Kronheimer-Mrowka's proof of the local Thom conjecture. \cite{MR1241873}

\begin{figure}
		\begin{subfigure}{1\textwidth}
			\centering
	\executeiffilenewer{product-sutured-manifold-product-disk.svg}{product-sutured-manifold-product-disk.pdf}{%
		\Inkscape -z -D --file="product-sutured-manifold-product-disk.svg" --export-pdf="product-sutured-manifold-product-disk.pdf" --export-latex}%
	\scalebox{1}{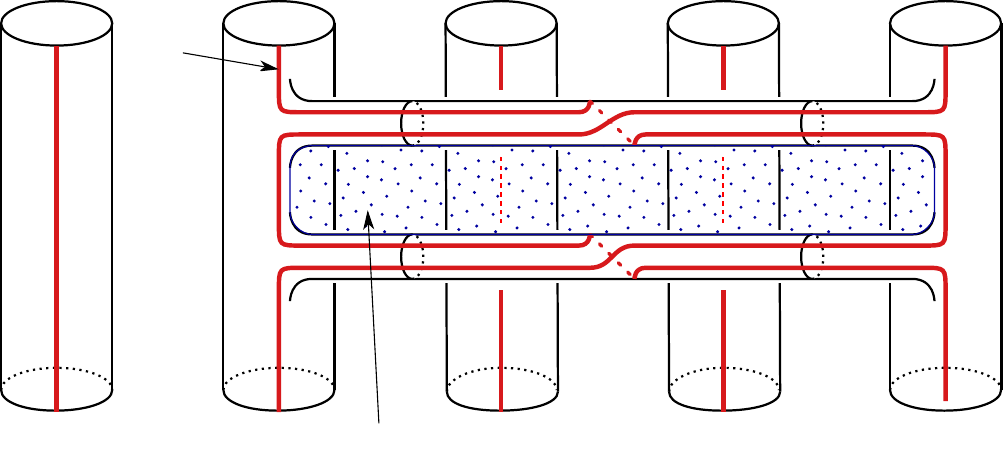}%

			\caption{The piece of $(Y_w, \gamma_w)$ corresponding to $a_{i,j}^2$. The dotted blue disk $D$ is the product disk and the solid red arcs are part of the sutures $\gamma_w$.}
			\label{fig:generatorsquaredcomplsuturedmfld}
		\end{subfigure}
		\par\bigskip
		\par\bigskip
		\par\bigskip
		\begin{subfigure}{1\textwidth}
			\centering
	\executeiffilenewer{product-sutured-manifold-product-disk-decomposed.svg}{product-sutured-manifold-product-disk-decomposed.pdf}{%
		\Inkscape -z -D --file="product-sutured-manifold-product-disk-decomposed.svg" --export-pdf="product-sutured-manifold-product-disk-decomposed.pdf" --export-latex}%
	\scalebox{.4}{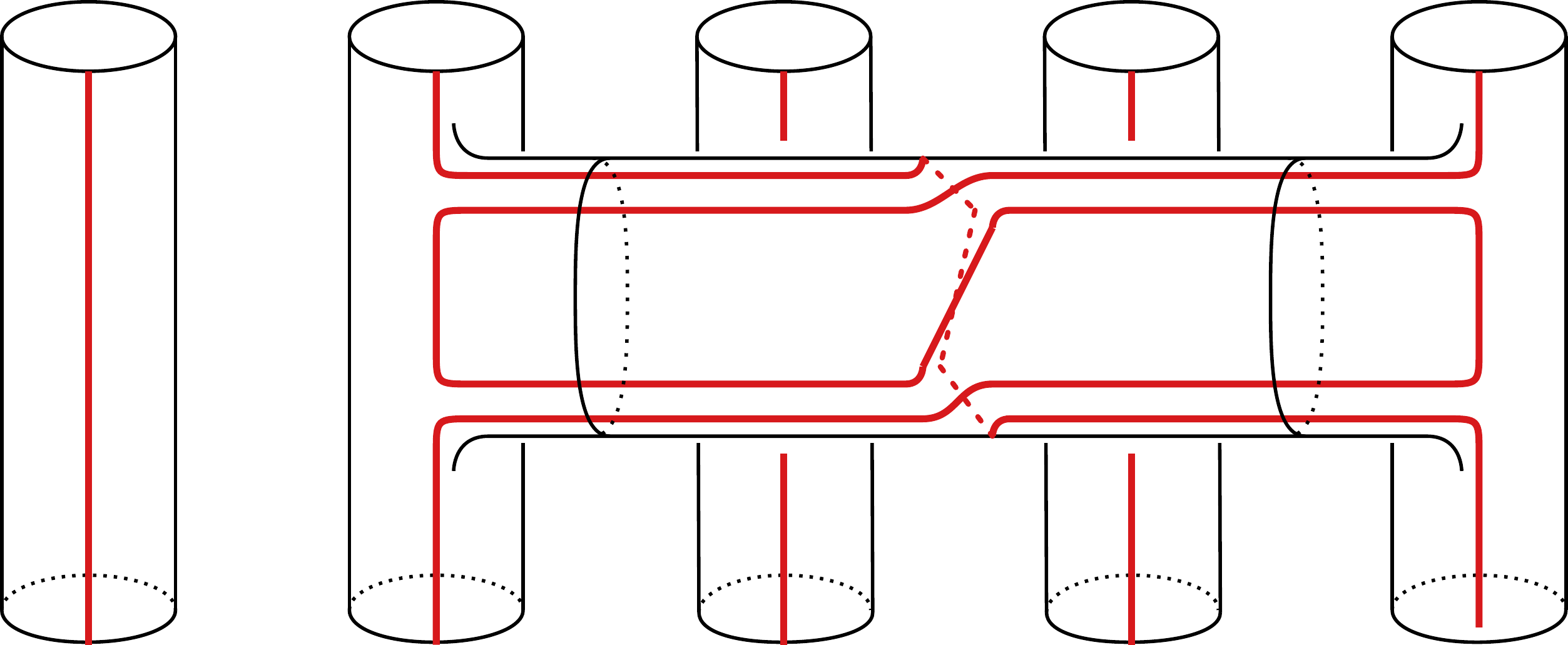}%

			\caption{Applying the product disk decomposition to $(Y_w, \gamma_w)$.}
			\label{fig:generatorsquaredcomplsuturedmflddecomposed}
		\end{subfigure}
		\par\bigskip
		\par\bigskip
		\par\bigskip
		\begin{subfigure}{1\textwidth}
			\centering
	\executeiffilenewer{product-sutured-manifold-product-disk-decomposed-simplified.svg}{product-sutured-manifold-product-disk-decomposed-simplified.pdf}{%
		\Inkscape -z -D --file="product-sutured-manifold-product-disk-decomposed-simplified.svg" --export-pdf="product-sutured-manifold-product-disk-decomposed-simplified.pdf" --export-latex}%
	\scalebox{.4}{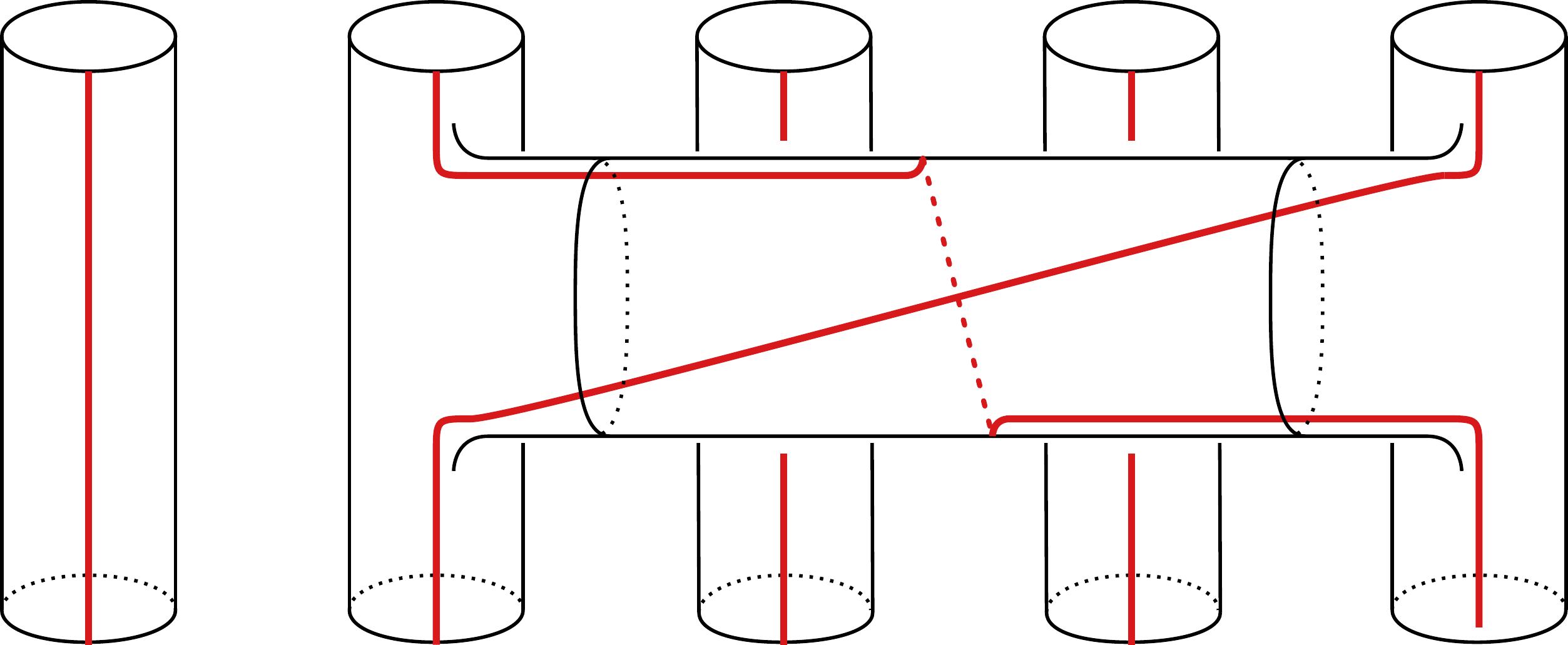}%

			\caption{An isotopy of the suture from Figure \ref{fig:generatorsquaredcomplsuturedmflddecomposed} shows that the decomposed sutured manifold is $(Y_{w'}, \gamma_{w'})$.}
			\label{fig:generatorsquaredcomplsuturedmfldsimplified}
		\end{subfigure}
		\caption{Proof of Lemma \ref{lemma:cancelgenerators}}
	\end{figure}

\begin{lemma}[Cancellation]
	\label{lemma:cancelgenerators}
	Suppose a braid word $w$ contains a term $a_{i,j}^2$, and let $w'$ be the word obtained from $w$ by replacing the $a_{i,j}^2$ with $a_{i,j}$ (i.e. a square of a generator is replaced with just the generator). Then the complementary sutured manifold $(Y_w, \gamma_w) = (S^3 \setminus (\Sigma_w \times I), \partial \Sigma_w)$ is a product if and only if $(Y_{w'}, \gamma_{w'}) = (S^3 \setminus (\Sigma_{w'} \times I), \partial \Sigma_{w'})$ is a product.
\end{lemma}
\begin{proof}
	According to Gabai (\cite{MR870705}, Lemma 2.2), if $(Y,\gamma) \rightsquigarrow (Y',\gamma')$ is a product disk decomposition, then $(Y,\gamma)$ is a product if and only if $(Y',\gamma')$ is a product. The lemma will follow immediately once we show that $(Y_{w'}, \gamma_{w'})$ is obtained from $(Y_w, \gamma_w)$ by a product disk decomposition.
	See figures \ref{fig:generatorsquaredcomplsuturedmfld}, \ref{fig:generatorsquaredcomplsuturedmflddecomposed}, \ref{fig:generatorsquaredcomplsuturedmfldsimplified} for the proof.
\end{proof}

Lemma \ref{lemma:cancelgenerators} and Theorem \ref{thm:maintheorem} were inspired by Ni's work on fibered 3-braids \cite{MR2597240}. The operation of canceling repeated generators of $B_3$ was called ``Untwisting" in Ni's work.

For strongly quasipositive braid closures the Bennequin surface is minimal genus and hence a fiber surface if the braid closure is fibered. This lemma then immediately implies that the closure of $\beta = w$ is fibered if and only if the closure of $\beta'$ is fibered, meaning we may cancel powers of generators in deciding if strongly quasipositive braid closures are fibered.

\section{Strongly quasipositive fibered braid closures}

\begin{lemma}
	Let $\beta = \delta P$, where $\delta = \sigma_{n-1} \sigma_{n-2} \dots \sigma_{1}$ is the Dual Garside element and $P$ is a BKL positive word. Then the braid closure $\hat{\beta}$ is fibered.
	\label{thm:mainlemma}
\end{lemma}
\begin{proof}
	We use induction on the word length $L = ln(P)$ of the BKL-positive word $P$ in the BKL generators. If $L = 0$, then $\beta = \delta^k$ is a non-split Artin-positive braid, whose closure is fibered by a result of Stallings \cite{MR520522}. 
	
	Otherwise, if $L > 0$, write $P = a_{s,t} P'$ for some generator $a_{s,t} \in S(P)$ and a BKL-positive word $P'$. Since $a_{s,t}$ is in the finishing set $F(\delta)$ of the dual Garside element $\delta$ by Fact \ref{fact:finishingsetofdualgarside}, we can write $\delta = P'' a_{s,t}$, where $P''$ is a BKL-positive word. Then 
	\begin{align*}
	\beta 
	&= \delta^k P \\
	&= \delta^{k-1} (P'' a_{r,s}) (a_{r,s} P') \\
	&= \delta^{k-1} (P'' (a_{r,s} a_{r,s}) P').
	\end{align*}
	Canceling the repeated generator $a_{r,s}$, we obtain the braid $\beta' = (\delta^{k-1} P'' a_{r,s} P') =\delta^k P'$. Note that $\ln(\beta') = L-1 < L$, and by induction we conclude that the braid closure $\hat{\beta}$ is fibered. By Fact \ref{fact:gabaifibredlinks} and Lemma \ref{lemma:cancelgenerators}, if $\hat{\beta'}$ is fibered, then $\hat{\beta}$ is also fibered.
\end{proof}

\begin{theorem}
	Let $\beta \in B_n$ be strongly quasipositive and let $\beta = \delta^k A_1 \dots A_m$ be its normal form. If $k \geq 1$, then the braid closure $\hat{\beta}$ is fibered.
	\label{thm:maintheorem}
\end{theorem}
\begin{proof}
	Apply Lemma \ref{thm:mainlemma} to $\beta = \delta P$ with $P = A_1 \dots A_m$.
\end{proof}
The converse of Theorem \ref{thm:maintheorem} is not true. In Section \ref{section:hopfplumbedbaskets} we explain a geometric classification of braid closures whose normal form contains a positive power of the Dual Garside element, and use this in Corollary \ref{cor:notallpositivefibredsqpcontaindelta} to show that the $(2,1)$ cable of the trefoil is a fibered strongly quasipositive braid closure whose normal form does not contain a positive power of the Dual Garside element. 

The following corollary is known for braids on 3 strands. \cite{MR2597240,2006math......6435S}

\begin{figure}
	\centering
	\executeiffilenewer{BKL-proof-of-cor-to-mainthm.svg}{BKL-proof-of-cor-to-mainthm.pdf}{%
		\Inkscape -z -D --file="BKL-proof-of-cor-to-mainthm.svg" --export-pdf="BKL-proof-of-cor-to-mainthm.pdf" --export-latex}%
	\scalebox{1}{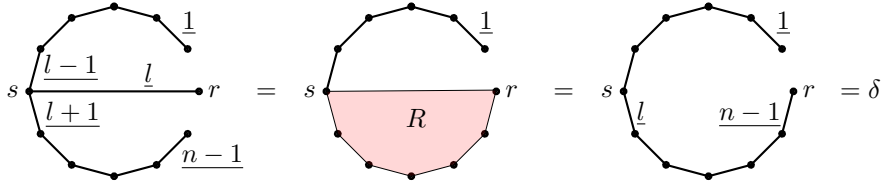}%

	\caption{We rewrite the braid $\beta_R$ corresponding to the polygon $R$ as $\beta_R = a_{s,s-1} \dots a_{r+1,r}$, where $l = s-r$ (mod $n$). }
	\label{fig:makingsqpfibredbyaddingcrossings}
\end{figure}

\begin{cor}
	After adding at most $n-2$ crossings to $\beta \in B_n$, every non-split strongly quasipositive braid closure $\hat{\beta}$ on $n$ strands becomes fibered.
	\label{lemma:makingsqpfibredbyaddingcrossings}
\end{cor}
\begin{proof}
	By the previous theorem, we may assume that the normal form is \\
	$\beta = A_1 \dots A_m$. Let $A_m = P a_{s,r}$ for a BKL-positive word $P$ and a generator $a_{s,r} \in F(A_m)$. Consider the strongly quasipositive braid $\beta'$ obtained from $\beta$ by adding $n-2$ crossings (all subscripts are mod $n$):
	\begin{align*}
	\beta' &= A_1 \dots A_{m-1} P (a_{r-2, r-1} a_{r-3, r-2} \dots a_{s+1, s}) a_{s,r} (a_{s,s-1} a_{s-1,s-2} \dots a_{r+2, r+1}) \\ 
	&=  A_1 \dots A_{m-1} P \delta \mbox{ ~ ~(see Figure \ref{fig:makingsqpfibredbyaddingcrossings})}
	\end{align*}
	The braid $\beta'$ is conjugate to $\delta A_1 \dots A_{m-1} P$ and hence its closure $\hat{\beta'}$ is fibered by Theorem \ref{thm:maintheorem}.
\end{proof}

Theorem \ref{thm:maintheorem} has the following probabilistic interpretation: Fixing the number of strands, the probability that a randomly generated strongly quasipositive braid word will contain $\delta$, or in fact any subword of fixed length, approaches 1 as the word length increases to $\infty$. This justifies our claim that almost all strongly quasipositive braid closures are fibered.

The above observations may lead the reader to conclude that \textit{all} strongly-quasipositive braids are fibered, but this is not true. In fact, on the class of strongly quasipositive links, all Seifert forms are realized, so the leading coefficient of the Alexander polynomial $\Delta_L$ may be arbitrary! \cite{MR2179266}

In particular, if $\Delta_L$ is not monic, then the link $L$ is not fibered. An example of a non-fibered strongly quasipositive braid closure is given by $\beta_3 = a_{3,1} a_{4,2} a_{3,1} a_{4,2}$. The Bennequin surface $\Sigma_{\beta_3}$ consists of two once-linked Hopf annuli. The complement of a thickening of $\Sigma_{\beta_3}$ is the Hopf-link exterior, which implies that the closure $\hat{\beta_3}$ is not fibered as the Hopf-link exterior is not a handlebody.

\section{Hopf-Plumbed baskets}
\label{section:hopfplumbedbaskets}

In \cite{MR1857666}, Rudolph constructs fibered links that arise as closures of certain homogeneous braids $\beta \in B_n$, generalizing earlier work of Stallings on closures of homogeneous braids in the Artin generators \cite{MR520522}. Rudolph constructs generating set $G(\mathcal{T})$ of $B_n$ as follows:

\begin{itemize}
	\item Let $\mathcal{T}$ be a tree with $n$ vertices and $n-1$ edges. Embed the tree into $\mathbb{C}$ with vertices at $1,2, \dots n \in \mathbb{R} \subset \mathbb{C}$ and edges in the lower-half plane. Note that the assumptions imply that every vertex $1, \dots, n$ is the endpoint of at least one edge. These trees are called espaliers.
	\item To an edge $e \in E(\mathcal{T})$ of $\mathcal{T}$ with endpoints at the vertices $r,s$ associate the BKL-generator $a(e) = a_{s,r} \in B_n$. 
	\item Let $G(\mathcal{T}) = \{ a(e)~|~e \in E(\mathcal{T}) \}$ be the set of $\mathcal{T}$ generators.
\end{itemize}

\begin{figure}
	\begin{subfigure}[t]{.49\textwidth}
		\centering
	\executeiffilenewer{Espalier_Arbitrary.svg}{Espalier_Arbitrary.pdf}{%
		\Inkscape -z -D --file="Espalier_Arbitrary.svg" --export-pdf="Espalier_Arbitrary.pdf" --export-latex}%
	\scalebox{1}{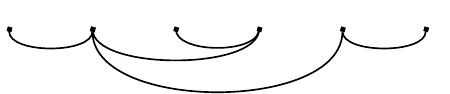}%

		\subcaption{$G(\mathcal{T}) = \{a_{2,1}, a_{4,2}, a_{5,2}, a_{4,3}, a_{6,5}\}$}
		\label{fig:espalier_arbitrary}
	\end{subfigure}
	\begin{subfigure}[t]{.49\textwidth}
		\centering
	\executeiffilenewer{Espalier_Max.svg}{Espalier_Max.pdf}{%
		\Inkscape -z -D --file="Espalier_Max.svg" --export-pdf="Espalier_Max.pdf" --export-latex}%
	\scalebox{1}{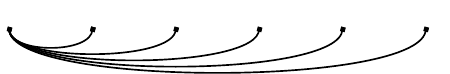}%

		\subcaption{$G(\mathcal{Y}_6) = \{a_{2,1}, a_{3,1}, a_{4,1}, a_{5,1}, a_{6,1}\}$}
		\label{fig:espalier_max}
	\end{subfigure}
	\caption{Espaliers and associated generating sets.}
	\label{fig:espaliers}
\end{figure}

See Figure \ref{fig:espaliers} for examples of espaliers and their generating sets. The espalier in Figure \ref{fig:espalier_max} is maximal \footnote{$\mathcal{Y}_n$ was denoted $\mathcal{Y}_X$ (for $X = \{1,2,\dots, n\}$) in \cite{MR1857666}. The definition of the partial order is given in Chapter 6 of Rudolph's paper \cite{MR1857666}.} among espaliers on $n$ vertices.

Let $H = \langle S | R\rangle$ be a group presentation with generating set $S$ and relations $R$. A word $w$ is \textbf{homogeneous} (\textbf{positive}) with respect to the generating set $S$ if every generator $s \in S$ occurs in $w$ with either only positive or only negative (only positive) powers. 
A word $\beta$ in the generating set $G(\mathcal{T})$ associated to a tree $\mathcal{T}$ is called a \textbf{$\mathcal{T}$-bandword}. A \textbf{strict} $\mathcal{T}$-bandword is a homogeneous $\mathcal{T}$-bandword $\beta$ such that every generator $g \in G(\mathcal{T})$ occurs in $\beta$.

\begin{definition}[\cite{MR1857666}]
	A surface $S$ is a Hopf-plumbed basket if $S$ is a plumbing of Hopf bands along arcs $\alpha_i \subset D$, all of which lie in $D$. A $(+)$ Hopf-plumbed basket is a Hopf-plumbed basket such that all plumbands are positive Hopf bands.
\end{definition}

\begin{theorem}
	A link $L$ is the boundary of $(+)$ Hopf plumbed basket if and only if $L$ admits a strongly quasipositive representative $\beta \in B_n$ which contains the Dual Garside element $\delta$.
	\label{thm:hopfplumbedbaskets}
\end{theorem}

In particular, this implies that the converse to Theorem \ref{thm:maintheorem} is not true. The author would like to thank Sebastian Baader for providing the following counterexample.

\begin{cor}
	There exist fibered strongly quasipositive braid closures which do not contain the Dual Garside element $\delta$.
	\label{cor:notallpositivefibredsqpcontaindelta}
\end{cor}
\begin{proof}
	Let $L$ be the $(2,1)$ cable of the trefoil, and let $F$ be the fiber surface for $L$. By work of Hedden \cite{MR2509749} and earlier work of Melvin-Morton \cite{MR859157}, $F$ is quasipositive but does not deplum a Hopf band. By definition, Hopf-plumbed baskets always deplum Hopf bands, so $F$ can not be a $(+)$ Hopf-plumbed basket. Theorem \ref{thm:hopfplumbedbaskets} then implies the boundary $L = \partial F$ does not admit a strongly quasipositive representative which contains the Dual Garside element $\delta$.
\end{proof}

\begin{lemma}
	Let 
	\[
		\beta = a_{r_1, 1} a_{r_2, 1} \dots a_{r_M, 1} \in B_n
	\] 
	be a positive $\mathcal{Y}_n$-bandword of word length $M$. Consider the sequence $(r_i) = (r_1, r_2, \dots, r_M)$. If there exists $1 \leq k < n$ and $1 \leq L \leq P \leq U \leq M$ such that
	\begin{itemize}
		\item $r_L = 1$, $r_P = k$, $r_U = n$,
		\item for all $i$ such that $1 \leq i \leq k$, there exists $Q(i)$ satisfying $L \leq Q(i) \leq P$ such that $r_{Q(i)} = i$ and further, $Q(i) < T \leq P$ implies that $r_{T} > r_{Q(i)}$,
		\item for all $i$ such that $k \leq i \leq n$, there exists $Q(i)$ satisfying $P \leq Q(i) \leq U$ such that $r_{Q(i)} = i$ and further, $P \leq T < Q(i)$ implies that $r_{T} < r_{Q(i)}$,
	\end{itemize}
	then $\beta$ contains the Dual Garside element $\delta$.
	\label{lemma:containsdelta}
\end{lemma}
\begin{proof}
	Informally, the idea of the proof is to slide $a_{i,1}$ ``down" for $i < k$, and slide  $a_{j,1}$ ``up" for $j > k$ so they form $\delta$. We introduce the following notation: Let $w(i,j)$ be the subword of $w$ from position $i$ to position $j$.\footnote{For example, take $z = a_{2,1} a_{6,3} a_{4,1} a_{4,2}$. Then $z(2,3) = a_{6,3} a_{4,1}$.} Now let $w$ be a positive $\mathcal{T}$-bandword for $\beta$ satisfying the statement of the lemma. We will prove that the subword $w(Q(1), Q(n)$ of $\beta$ contains $\delta$.
			For convenience, define $Q(k) = P$. For $s \in S = \{r_{Q(i)}, \dots, r_{Q(i+1) - 1 } \}$, the  relations in the braid group imply that $a_{i,1} a_{s,1} = a_{i, s} a_{i,1}$ as $s > i$ by our assumptions. Consider
			
			\begin{align*}
				a_{i,1} w(Q(i), Q(i+1) - 1) &= a_{i,1} a_{r_{Q(i)},1} a_{r_{Q(i)+1},1}  \dots a_{r_{Q(i+1)-1}, 1} \\
				&= \underbrace{a_{r_{Q(i)},i} a_{r_{Q(i)+1},i} \dots a_{r_{Q(i+1)}-1, i}}_{w'_i} a_{i,1},
			\end{align*}
			and note that $w'_i$ commutes with $a_{j,1}$ for $j < i$. We apply this equation in the next step:
			
			\begin{align*}
				w(Q(1), Q(k) - 1)
				&= \prod_i^{k-1} a_{i,1} w(Q(i), w(Q(i+1) - 1) \\
				&= \prod_i^{k-1} w'_{i} a_{i,1} \\
				&= w'_1 w'_2 \dots w'_{k-1} a_{1,1} a_{2,1} \dots a_{k-1,1} \\
				&= T_L a_{1,1} a_{2,1} \dots a_{k-1,1},
			\end{align*}
			for a BKL-positive word $T_L$. Moreover, a similar reasoning shows that for some BKL-positive word $T_U$,
			\[
				w(Q(k)+1, Q(n)) = a_{k+1,1} a_{k+1,1} \dots a_{n,1} T_U.
			\]
		To finish the proof, we note that the dual Garside element $\delta$ can be written as $\delta = a_{1,1 }a_{2,1} a_{3,1} \dots a_{n,1}$. We can now show that the subword $w(Q(1), Q(n)$ of $\beta$ contains $\delta$:
		\begin{align*}
			w(Q(1), Q(n)) &= w(Q(1), Q(k) - 1) Q(k) w(Q(k)+1, Q(n)) \\
			&=  T_L a_{1,1} a_{2,1} \dots a_{k-1,1} a_{k,1} a_{k+1,1} \dots a_{n,1} T_U \\
			&= T_L \delta T_U.
		\end{align*}
	\textbf{Remark} The converse is true: If $\beta$ contains the Dual Garside element $\delta$ then there is a word $w$ for $\beta$ which satisfies the assumptions of the lemma.
\end{proof}

\begin{figure}
	\begin{subfigure}[t]{.24\textwidth}
		\centering
	\executeiffilenewer{ProofThatAHopfPlumbedBasketContainsDelta_Step1.svg}{ProofThatAHopfPlumbedBasketContainsDelta_Step1.pdf}{%
		\Inkscape -z -D --file="ProofThatAHopfPlumbedBasketContainsDelta_Step1.svg" --export-pdf="ProofThatAHopfPlumbedBasketContainsDelta_Step1.pdf" --export-latex}%
	\scalebox{0.5}{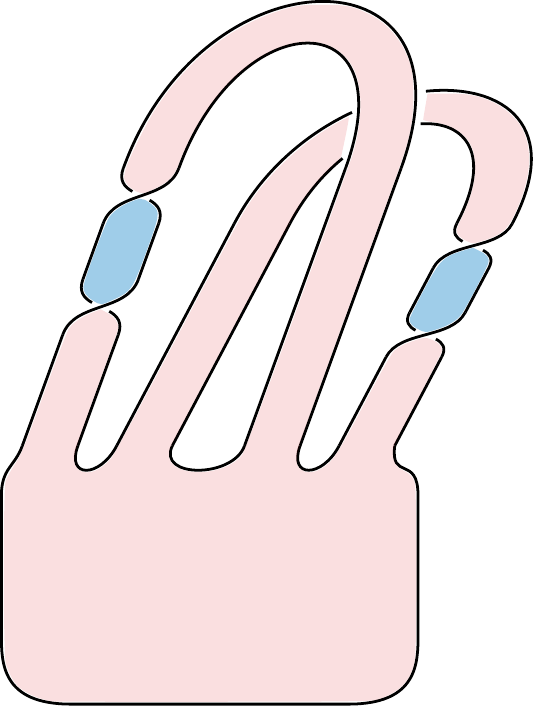}%

		\subcaption{$F$}
		\label{fig:hopfplumbedbasket}
	\end{subfigure}
	\begin{subfigure}[t]{.23\textwidth}
		\centering
	\executeiffilenewer{ProofThatAHopfPlumbedBasketContainsDelta_Step2.svg}{ProofThatAHopfPlumbedBasketContainsDelta_Step2.pdf}{%
		\Inkscape -z -D --file="ProofThatAHopfPlumbedBasketContainsDelta_Step2.svg" --export-pdf="ProofThatAHopfPlumbedBasketContainsDelta_Step2.pdf" --export-latex}%
	\scalebox{0.6}{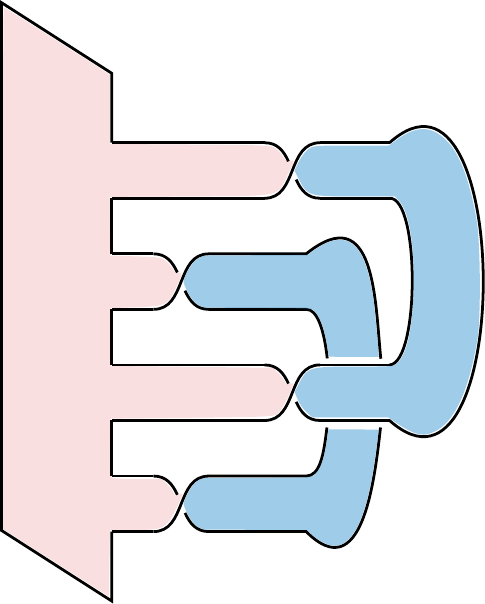}%

		\subcaption{Isotopy of $F$}
		\label{fig:hopfplumbedbasket_isotoped}
	\end{subfigure}
	\begin{subfigure}[t]{.26\textwidth}
		\centering
	\executeiffilenewer{ProofThatAHopfPlumbedBasketContainsDelta_Step3.svg}{ProofThatAHopfPlumbedBasketContainsDelta_Step3.pdf}{%
		\Inkscape -z -D --file="ProofThatAHopfPlumbedBasketContainsDelta_Step3.svg" --export-pdf="ProofThatAHopfPlumbedBasketContainsDelta_Step3.pdf" --export-latex}%
	\scalebox{0.6}{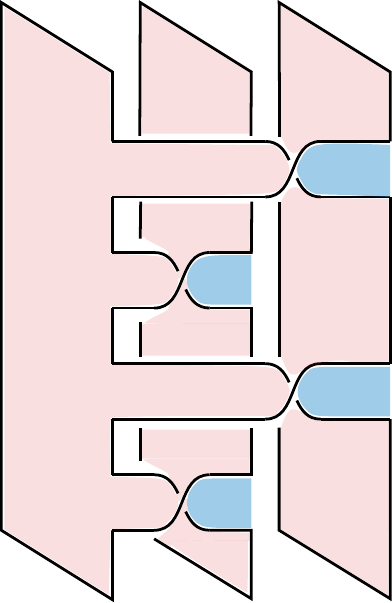}%

		\subcaption{Isotopy to a quasipositive surface}
		\label{fig:hopfplumbedbasket_isotopedtosqpbraid}
	\end{subfigure}
	\begin{subfigure}[t]{.24\textwidth}
		\centering
	\executeiffilenewer{ProofThatAHopfPlumbedBasketContainsDelta_Step4.svg}{ProofThatAHopfPlumbedBasketContainsDelta_Step4.pdf}{%
		\Inkscape -z -D --file="ProofThatAHopfPlumbedBasketContainsDelta_Step4.svg" --export-pdf="ProofThatAHopfPlumbedBasketContainsDelta_Step4.pdf" --export-latex}%
	\scalebox{1}{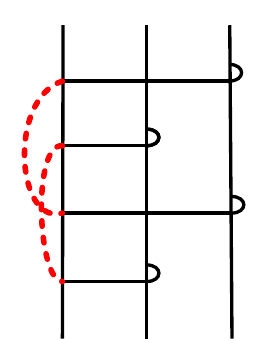}%

		\subcaption{Charged fence diagram}
		\label{fig:hopfplumbedbasket_chargedfencediagram}
	\end{subfigure}
	\caption{The fiber surface $F$ for the trefoil and a charged fence diagram for $F$. The positive Hopf band $A_+$ plumbed to $\alpha_i \in D$ is represented by the strand $i > 1$ together with the two horizontal lines connecting strand $1$ to $i$ and the dotted arc $\alpha_i$.}
	\label{fig:hopfplumbedbasketfortrefoil}
\end{figure}

For the proof of Theorem \ref{thm:hopfplumbedbaskets}, it will be convenient to use so-called charged fence diagrams, which represent quasipositive surfaces $S$ and define a braid word for the strongly quasipositive link $\partial S$. See figures \ref{fig:chargedfencediagram} and \ref{fig:hopfplumbedbasket_chargedfencediagram} for examples of charged fence diagrams, and \cite{MR1452826} for a thorough discussion of quasipositive surfaces and charged fence diagrams. 

\begin{figure}
	\begin{subfigure}[t]{.48\textwidth}
		\centering
	\executeiffilenewer{PlumbingHopfBandsToDeltaVersion2_Step1.svg}{PlumbingHopfBandsToDeltaVersion2_Step1.pdf}{%
		\Inkscape -z -D --file="PlumbingHopfBandsToDeltaVersion2_Step1.svg" --export-pdf="PlumbingHopfBandsToDeltaVersion2_Step1.pdf" --export-latex}%
	\scalebox{0.8}{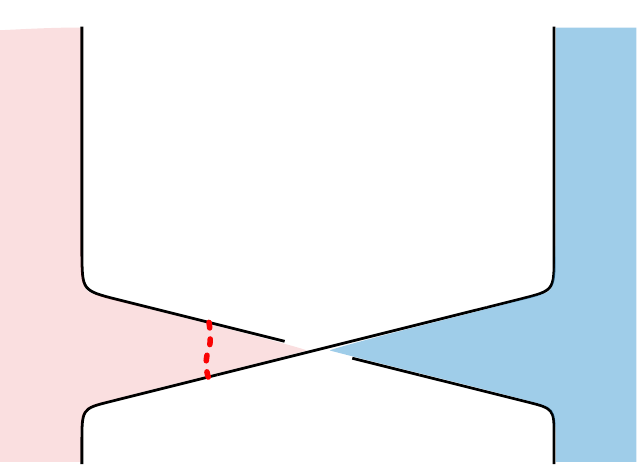}%

		\subcaption{A piece of the disk $D$ and $\alpha \subset D$.}
		\label{fig:attachinghopfbandtodeltastep1}
	\end{subfigure}
	\begin{subfigure}[t]{.48\textwidth}
		\centering
	\executeiffilenewer{PlumbingHopfBandsToDeltaVersion2_Step2.svg}{PlumbingHopfBandsToDeltaVersion2_Step2.pdf}{%
		\Inkscape -z -D --file="PlumbingHopfBandsToDeltaVersion2_Step2.svg" --export-pdf="PlumbingHopfBandsToDeltaVersion2_Step2.pdf" --export-latex}%
	\scalebox{0.8}{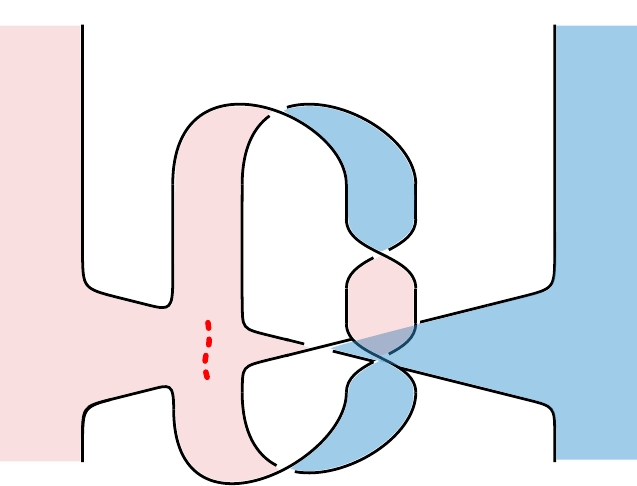}%

		\subcaption{Plumbing a positive Hopf along $\alpha$.}
		\label{fig:attachinghopfbandtodeltastep2}
	\end{subfigure}
	\vskip\baselineskip
	\begin{subfigure}[t]{.48\textwidth}
		\centering
	\executeiffilenewer{PlumbingHopfBandsToDeltaVersion2_Step3.svg}{PlumbingHopfBandsToDeltaVersion2_Step3.pdf}{%
		\Inkscape -z -D --file="PlumbingHopfBandsToDeltaVersion2_Step3.svg" --export-pdf="PlumbingHopfBandsToDeltaVersion2_Step3.pdf" --export-latex}%
	\scalebox{0.8}{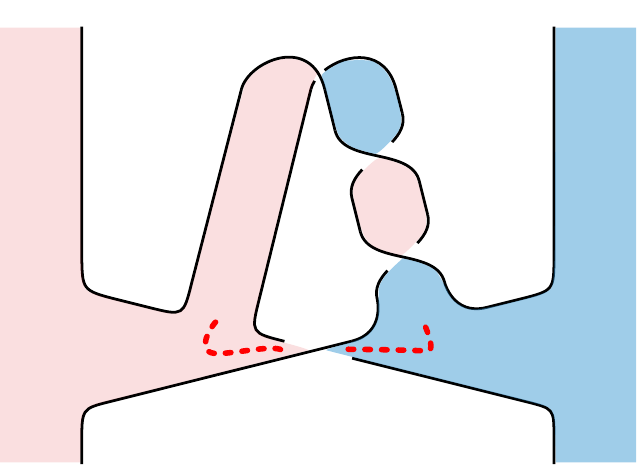}%

		\subcaption{Isotopy.}
		\label{fig:attachinghopfbandtodeltastep3}
	\end{subfigure}
	\begin{subfigure}[t]{.48\textwidth}
		\centering
	\executeiffilenewer{PlumbingHopfBandsToDeltaVersion2_Step4.svg}{PlumbingHopfBandsToDeltaVersion2_Step4.pdf}{%
		\Inkscape -z -D --file="PlumbingHopfBandsToDeltaVersion2_Step4.svg" --export-pdf="PlumbingHopfBandsToDeltaVersion2_Step4.pdf" --export-latex}%
	\scalebox{0.8}{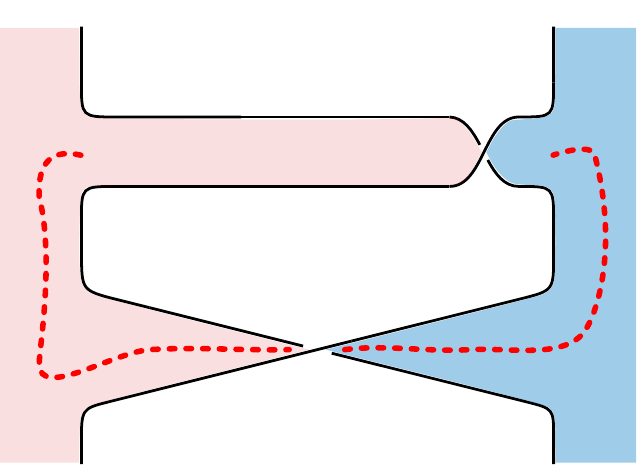}%

		\subcaption{Further isotopy.}
		\label{fig:attachinghopfbandtodeltastep4}
	\end{subfigure}
	\caption{Plumbing a positive Hopf band along $\alpha \in D$ ``adds the generator $a_{s,r}$."}
	\label{fig:proofplumbhopfbandtodelta}
\end{figure}

\begin{lemma}
	Let $\beta = a_{s_1,r_1} \dots a_{s_k,r_k} \delta \in B_n$ be a braid which contains the Dual Garside element $\delta$. Then the fiber surface for $\hat{\beta}$ is a $(+)$ Hopf-plumbed basket.
	\label{lemma:dualgarsideimplieshopfplumbedbasket}
\end{lemma}
\begin{proof}
	We construct a $(+)$ Hopf-plumbed basket whose boundary is $\hat{\beta}$. Let
	\[
		D = (D_1 \cup \dots \cup D_n) \cup (T_{\sigma_1} \cup \dots \cup T_{\sigma_{n-1}}),
	\] be the Bennequin surface for the subword $\delta = \sigma_{n-1} \dots \sigma_1$ of $\beta$. Write $\delta = a_{s_1,r_1} P$ for a BKL-positive word $P$ and isotope $D$ accordingly. 
	
	Consider the surface $D \star_{\alpha_1} A_+$ given by plumbing a positive Hopf band $A_+$ to this disk $D$ along the arc $\alpha_1 = \alpha$ in Figure \ref{fig:attachinghopfbandtodeltastep1}. The surface $D \star_{\alpha_1} A_+$ is pictured in Figure \ref{fig:attachinghopfbandtodeltastep2}, and the sequence of isotopies in figures \ref{fig:attachinghopfbandtodeltastep3} and \ref{fig:attachinghopfbandtodeltastep4} show that $D \star_{\alpha_1} A_+$ is the fiber surface for the closure of the braid $a_{s_1, r_1} \delta$. Informally, this plumbing ``adds the generator $a_{s,r}$ to the braid word."

	Now repeat this process by plumbing positive Hopf bands for the generators $(a_{s_2,r_2}), \dots , (a_{s_k,r_k})$. By construction, the fiber surface $((D \star_{\alpha_1} A_+)\dots \star_{\alpha_k} A_+)$ is a $(+)$ Hopf-plumbed basket bounding $\hat{\beta}$.
\end{proof}

\begin{figure}
	\begin{subfigure}[t]{.30\textwidth}
		\centering
	\executeiffilenewer{PlumbingHopfBandsToDelta_Example.svg}{PlumbingHopfBandsToDelta_Example.pdf}{%
		\Inkscape -z -D --file="PlumbingHopfBandsToDelta_Example.svg" --export-pdf="PlumbingHopfBandsToDelta_Example.pdf" --export-latex}%
	\scalebox{1}{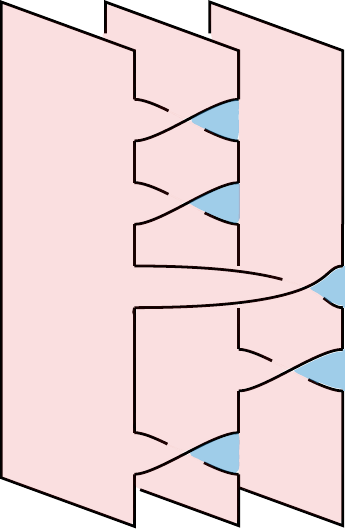}%

		\subcaption{Fiber surface $F$ bounding $\beta = a_{2,1} a_{2,1} a_{3,1} \delta$.}
		\label{fig:containsDeltaIsABasketStep1}
	\end{subfigure}
	\begin{subfigure}[t]{.36\textwidth}
		\centering
	\executeiffilenewer{PlumbingHopfBandsToDelta_ExampleStep2.svg}{PlumbingHopfBandsToDelta_ExampleStep2.pdf}{%
		\Inkscape -z -D --file="PlumbingHopfBandsToDelta_ExampleStep2.svg" --export-pdf="PlumbingHopfBandsToDelta_ExampleStep2.pdf" --export-latex}%
	\scalebox{1}{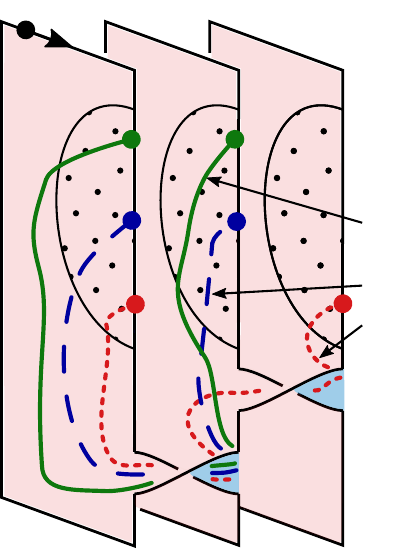}%

		\subcaption{$D = D_1 \cup D_2 \cup D_3 \cup T_1 \cup T_2$ and the plumbing arcs $\alpha_i \subset D$. $F = ((D \star_{\alpha_1} A_+) \star_{\alpha_2} A_+) \star_{\alpha_3} A_+$.}
		\label{fig:containsDeltaIsABasketStep2}
	\end{subfigure}
	\begin{subfigure}[t]{.32\textwidth}
		\centering
	\executeiffilenewer{UntwistedDelta.svg}{UntwistedDelta.pdf}{%
		\Inkscape -z -D --file="UntwistedDelta.svg" --export-pdf="UntwistedDelta.pdf" --export-latex}%
	\scalebox{1}{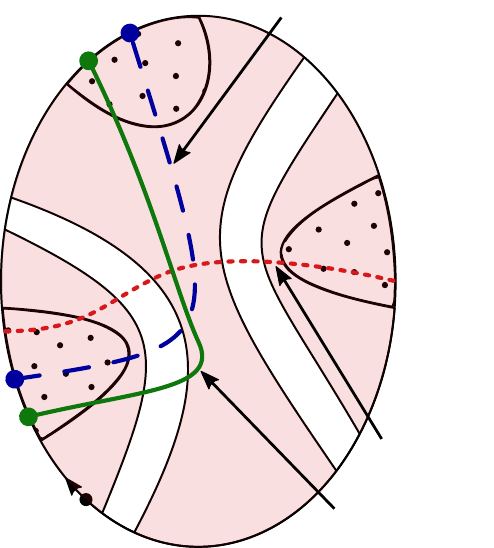}%

		\subcaption{Folding structure for $D$.}
		\label{fig:containsDeltaIsABasketStep3}
	\end{subfigure}
	\caption{}
	\label{fig:containsDeltaIsABasket}
\end{figure}

Given a braid $\beta = a_{r_1,s_1} \dots a_{r_k,s_k} \delta \in B_n$ representing $L = \hat{\beta}$, the lemma provides an explicit construction of a $(+)$ Hopf-plumbed basket $(D \star_{\alpha_1} A_+) \dots \star_{\alpha_k} A_+$ whose boundary is the link $L$. The normal form of the representative $\beta \in B_n$ induces additional structure on $D$, which is pictured in Figure \ref{fig:containsDeltaIsABasket}.
\begin{itemize}
	\item The disk $D$ is partitioned into disjoint disks $D_i$ for $1 \leq i \leq n$ and $T_j$ for $1 \leq j < n$ as in Figure \ref{fig:containsDeltaIsABasketStep3}. The disks $D_1$ and $D_n$ are distinguished by noting that $\partial D \cap D_i$ has a single component for $i = 1$ and $i=n$ but two components for $1 < i < n$. There are two components in $\partial D \cap T_j$ for every $i \leq j < n$.
	
	\item There exist disks $AT_i \subset D_i$ which satisfy
	(1) $\cup \partial \alpha_j \subset \cup \partial AT_i$ and (2) $AT_i \cap \partial D_i$ is connected, where $\{\alpha_j\}$ is an ordered set of compatible arcs (defined below).  We refer to $AT_i$ as an ``attaching region." In figures \ref{fig:containsDeltaIsABasketStep2} and \ref{fig:containsDeltaIsABasketStep3}, the attaching regions $AT_i$ are the dotted regions contained in $D_i$.
	
	\item Orient the arcs $\gamma_i = AT_i \cap \partial D_i$ as submanifolds of $\partial D$ endowed with the boundary orientation. We say that an ordered set $\{ \tau_j \subset D\}$ of properly embedded arcs whose endpoints are contained in $\cup_i \gamma_i$ is \textbf{compatible} if (1) $|\gamma_i \cap \partial \tau_j | < 2$ for all $i,j$; and (2) for $k > l$, if $a \in \partial \tau_k \cap \gamma_i$ and $b \in \partial \tau_l \cap \gamma_i$, then $a < b$ in the order on $\gamma_i \simeq (0,1)$ induced by the orientation.
\end{itemize}

The construction in the proof of the lemma implies that the plumbing arcs $\{ \alpha_i \}$ are compatible and that the disks $D_i$ ($T_j$, respectively) correspond to the disks $D_i$ (twisted bands for $\sigma_j= a_{j+1,j}$, respectively) in the Bennequin surface $\Sigma_\delta$.

We call such a structure a \textbf{folding structure} and denote it by $(D, D_i, \alpha_j)$, where $D_i$ are disks in a partition of a disk $D$ as described above and $\{\alpha_j\}$ is an ordered set of compatible arcs. To justify the omission of $T_j$, $AT_i$ and $\gamma_i$ from the notation, we note that a partition satisfying the conditions is completely determined by the disks $D_i$, essentially by the smooth Jordan curve theorem. The attaching regions $AT_i$ and the arcs $\gamma_i \subset \partial AT_i$ are defined in terms of the disks $D_i$ and the ordered set of compatible arcs. A priori, $(D, D_i, \alpha_j)$ may not be arising from a braid $\beta \in B_n$ whose normal form contains $\delta$ through the construction described in the proof of Lemma \ref{lemma:dualgarsideimplieshopfplumbedbasket}, however, the following lemma asserts this is the case and justifies the term ``folding structure."

\begin{lemma}
	Let $(D, D_i, \alpha_j)$ be a folding structure. Then there exists a braid $\beta \in B_n$ which contains the Dual Garside $\delta$ such that the folding structure is induced by $\beta$.
	\label{lemma:foldingstructurerealized}
\end{lemma}
\begin{proof}
	Let $AT_i$ be the attaching regions for the folding structure and let $k = |\{\alpha_j\}|$ be the number of arcs and $n = |\{D_i\}|$ be the number of disks. Define $r(j)$ and $s(j)$ by $\partial \alpha_j \in AT_{r(j)} \cup AT_{s(j)}$. Let $\beta = a_{r(1), s(1)} \dots a_{r(k), s(k)} \delta \in B_n$. It is straightforward to verify that the folding structure induced by $\beta$ is $(D, D_i, \alpha_j)$.
\end{proof}

\begin{lemma}
	The boundary $\partial F$ of a $(+)$ Hopf plumbed basket $F$ admits a representative which contains the Dual Garside element $\delta$.
	\label{lemma:turninghopfplumbedbasketintochargedfencediagram}
\end{lemma}

\begin{figure}
	\begin{subfigure}[t]{.33\textwidth}
		\centering
	\executeiffilenewer{BasketToDelta_Basket.svg}{BasketToDelta_Basket.pdf}{%
		\Inkscape -z -D --file="BasketToDelta_Basket.svg" --export-pdf="BasketToDelta_Basket.pdf" --export-latex}%
	\scalebox{1}{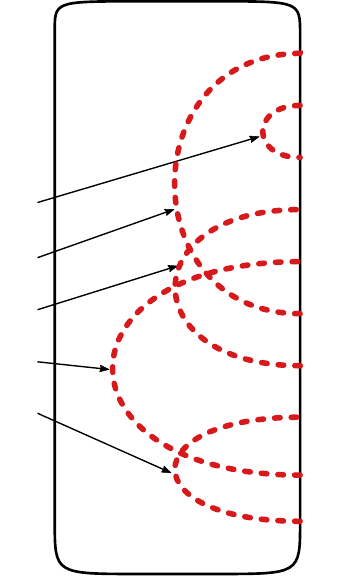}%

		\subcaption{Plumbing arcs $\{\alpha_j \subset D\}$ for $F$.}
		\label{fig:BasketContainsDeltaExample1}
	\end{subfigure}
	\begin{subfigure}[t]{.33\textwidth}
		\centering
	\executeiffilenewer{BasketToDelta_FoldingStructure.svg}{BasketToDelta_FoldingStructure.pdf}{%
		\Inkscape -z -D --file="BasketToDelta_FoldingStructure.svg" --export-pdf="BasketToDelta_FoldingStructure.pdf" --export-latex}%
	\scalebox{1}{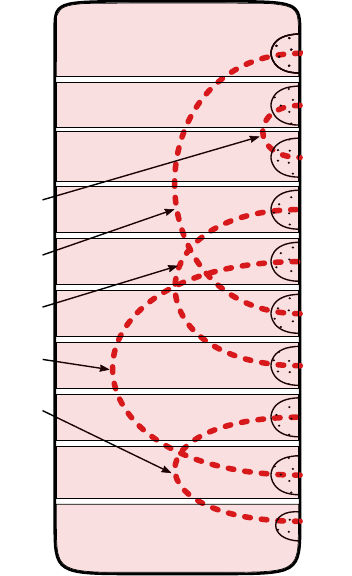}%

		\subcaption{Folding structure.}
		\label{fig:BasketContainsDeltaExample2}
	\end{subfigure}
	\begin{subfigure}[t]{.32\textwidth}
		\centering
	\executeiffilenewer{BasketToDelta_ChargedFenceDiagram.svg}{BasketToDelta_ChargedFenceDiagram.pdf}{%
		\Inkscape -z -D --file="BasketToDelta_ChargedFenceDiagram.svg" --export-pdf="BasketToDelta_ChargedFenceDiagram.pdf" --export-latex}%
	\scalebox{1}{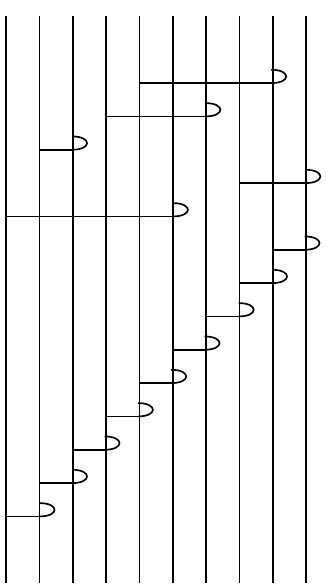}%

		\subcaption{The braid $\hat{\beta} = \partial F$ given by Lemma \ref{lemma:foldingstructurerealized}.}
		\label{fig:BasketContainsDeltaExample3}
	\end{subfigure}
	\caption{Proof of Lemma \ref{lemma:turninghopfplumbedbasketintochargedfencediagram}. A Hopf-plumbed basket $F$ is described by an ordered collection of arcs $\{\alpha_j \subset D\}$. The attaching regions (dotted) for the folding structure in \ref{fig:BasketContainsDeltaExample2} contain a single endpoint each to make $\{\alpha_j\}$ compatible.}
	\label{fig:BasketContainsDelta}
\end{figure}

\begin{proof}
	The idea of the proof is to define a folding structure $(D, D_i, \alpha_j)$ and appeal to the previous lemma. We begin by noting that according to Section 3 of \cite{MR1857666}, a Hopf-plumbed basket is completely described by the plumbing arcs $S = \{\alpha_j \subset D \}$ and the order in which the Hopf bands are plumbed. By re-indexing $S$, we may assume the Hopf bands are plumbed in increasing order of the plumbing arcs' indices.
	
	Let $k = | S |$ be the number of arcs and let $p_1, p_2, \dots, p_{2k} \in \partial D$ be the endpoints of the arcs $\{ \alpha_j\}$ in some cyclic order. Let $\tau \subset \partial D$ be the arc with endpoints at $p_1$ and $p_{2k}$ which is disjoint from all other endpoints $p_i$. Choose points $q_2, \dots q_{2k-1} \in \tau$ in increasing order on $\tau \simeq (0,1)$.
	
	We now define the disks $D_i$ of the folding structure. For $i=1$ and $i=2k$, let $D_i$ be a small disk containing $p_i$; and for $1 < i < 2k$, let $D_i$ be a tubular neighborhood of an arc connecting $p_i$ to $q_i$. This is illustrated in Figure \ref{fig:BasketContainsDeltaExample2}. The arcs $\{\alpha_j\}$ are compatible, as each attaching region contains exactly one endpoint, which implies that the compatibility condition is vacuous. By Lemma \ref{lemma:foldingstructurerealized}, the folding structure $(D, D_i, \alpha_j)$ is induced by a braid $\beta \in B_{2k}$ which contains the Dual Garside element $\delta$.
\end{proof}

\begin{cor}
	The boundary $\partial F$ of a $(+)$ Hopf-plumbed basket $F$ with plumbing arcs $\{\alpha_i\}$ is the closure of a strongly quasipositive braid $\beta = \gamma \delta$, where $\gamma$ is an unlink on $k = |\{\alpha_i\}|$ components.
\end{cor}

\begin{proof}[Proof of Theorem \ref{thm:hopfplumbedbaskets}]
	 Immediate from Lemma \ref{lemma:dualgarsideimplieshopfplumbedbasket} and Lemma \ref{lemma:turninghopfplumbedbasketintochargedfencediagram}.
\end{proof}

\begin{theorem}[\cite{MR1857666}, Theorems 6.1.6, 6.2.4]
	If $S$ is a Hopf-plumbed basket, then there is an espalier $\mathcal{T}$ and a strict homogeneous $\mathcal{T}$-bandword $\beta$ such that $\hat{\beta} = \partial S$. Conversely, the fiber surface of a strict homogeneous $\mathcal{T}$-bandword is a Hopf-plumbed basket.
	\label{theorem:Rudolphmain}
\end{theorem}

See Figure \ref{fig:hopfplumbedbasketfortrefoil} for the idea of the proof in case of $\mathcal{Y}_n$-bandwords. It is straightforward to check that Rudolph's proof holds when one restricts to strict positive $\mathcal{T}$-bandwords, which are strongly quasipositive braids:

\begin{theorem}[\cite{MR1857666}]
	If $S$ is a $(+)$ Hopf-plumbed basket, then there is an espalier $\mathcal{T}$ and a positive strict $\mathcal{T}$-bandword $\beta$ such that $\hat{\beta} = \partial S$. Conversely, the fiber surface of a positive strict $\mathcal{T}$-bandword is a $(+)$ Hopf-plumbed basket.
	\label{theorem:Rudolphmainsqp}
\end{theorem}

\begin{cor}
	A braid closure $\hat{\beta}$ can be represented by a positive strict $\mathcal{T}$-bandword if and only if it can be represented by a strongly quasipositive braid which contains the Dual Garside element.
	\label{cor:deltaequalsbandword}
\end{cor}
\begin{proof}
	Immediate from Theorem \ref{thm:hopfplumbedbaskets} and Theorem \ref{theorem:Rudolphmainsqp}.
\end{proof}

\begin{theorem}
	A non-split positive braid closure $L$ can be represented by a strongly quasipositive braid which contains a positive power of the Dual Garside element $\delta$.
	\label{thm:corpositivebraids}
\end{theorem}

\begin{proof}
	Let $\mathcal{T}_p$ be the espalier with $n$ vertices and with edges connecting $i$ to $i+1$ for $1 \leq i < n$. The set of generators $G(\mathcal{T}_p) = \{ a_{2,1} = \sigma_1, \dots, a_{n,n-1} = \sigma_{n-1}\}$ is the set of Artin generators of $B_n$, which shows that a non-split positive braid is a positive strict $\mathcal{T}_p$-bandword. Now apply Corollary \ref{cor:deltaequalsbandword}.
\end{proof}

\bibliographystyle{halpha}
\bibliography{fibered_sqp_braids}

\end{document}